\documentclass{gtart}     

\input gtmonout
\volumenumber{2}
\volumeyear{1999}
\volumename{Proceedings of the Kirbyfest}
\pagenumbers{1}{21}
\papernumber{1}
\received{24 September 1998}\revised{14 November 1999}
\published{22 November 1999}

\newtheorem{thm}{Theorem}[section]    
\newtheorem{cor}[thm]{Corollary}  

\theoremstyle{definition}
\newtheorem{defn}[thm]{Definition}    
\theoremstyle{definition}
\newtheorem{exam}[thm]{Example}
\newtheorem*{rem}{Remark}             
\newtheorem*{rems}{Remarks}        
\theoremstyle{definition}
\newtheorem*{nota}{Notation}             

\begin{document}
\title{Combinatorial Dehn surgery on cubed\\and Haken $3$--manifolds}
\asciititle{Combinatorial Dehn surgery on cubed and Haken 3-manifolds}
\shorttitle{Combinatorial Dehn surgery on cubed and Haken 3--manifolds}

\authors{I\thinspace R Aitchison\\J\thinspace H Rubinstein}
\asciiauthors{IR Aitchison and JH Rubinstein}

\address{Department of Mathematics and Statistics, University of Melbourne\\
Parkville, Vic 3052, Australia}
\email{iain@ms.unimelb.edu.au, rubin@ms.unimelb.edu.au}

\begin{abstract}  
A combinatorial condition is obtained for when immersed or embedded
incompressible surfaces in compact
3--manifolds with tori boundary components remain incompressible
after Dehn surgery. A combinatorial characterisation of hierarchies is
described. A new
proof is given of the topological rigidity theorem of Hass and Scott for
3--manifolds
containing immersed incompressible surfaces, as found in cubings of
non-positive curvature.
\end{abstract}
\asciiabstract{A combinatorial condition is obtained for when immersed 
or embedded incompressible surfaces in compact 3--manifolds with tori
boundary components remain incompressible after Dehn surgery. A
combinatorial characterisation of hierarchies is described. A new
proof is given of the topological rigidity theorem of Hass and Scott
for 3--manifolds containing immersed incompressible surfaces, as found
in cubings of non-positive curvature.}

\maketitle

\section{Introduction}\label{S:Intro}

 An important property of the class of $3$--dimensional manifolds with
cubings of
non-positive curvature is that they contain `canonical' immersed
incompressible
surfaces (cf \cite{AR1}). In particular these surfaces satisfy the
$4$--plane and
$1$--line
properties of Hass and Scott (cf \cite{HS}) and so any $P^2$--irreducible
$3$--manifold which
is homotopy equivalent to such a $3$--manifold is homeomorphic to it
(topological rigidity).
In this paper we study the behaviour of these canonical surfaces under Dehn
surgery on
 knots and links. A major objective here is to show that in the case of a
cubed manifold
with tori as boundary components, there is a simple criterion to tell if a
canonical surface
 remains incompressible after a particular Dehn surgery. This result is
very much
 in the spirit of the negatively curved Dehn surgery of Gromov and Thurston
(cf \cite{BH})
 and was announced in \cite{AR2}.  Many examples are given in \cite{AR3}.

The key lemma required is a combinatorial version of Dehn's lemma and the
loop theorem for immersed surfaces of the type considered by Hass and Scott
 with an extra condition --- the triple point property. We are able to give a
 simplified proof of the rigidity theorem of Hass and Scott for $3$--manifolds
 containing immersed incompressible surfaces with this additional condition.

 By analogy with the combinatorial Dehn's lemma and the loop theorem, we are
 also able to find a combinatorial characterisation of the crucial idea of
 hierarchies in $3$--manifolds, as used by Waldhausen in his solution to the
word
problem \cite{Wa1}. In particular, if a list of embedded surfaces is given
 in a $3$--manifold, with boundaries on the previous surfaces in the list, then
a simple condition determines whether all the surfaces are incompressible and
mutually boundary incompressible. This idea can be used to determine if
 the hierarchy persists after Dehn surgery --- for example, if there are
 several cusps in a $3$--manifold then we can tell if Dehn surgery on all
 but one cusp makes the final cusp remain incompressible (cf
\cite{AR2}).

One should also note the result of Mosher \cite{Mo} that the fundamental
 group of an atoroidal cubed $3$--manifold is word hyperbolic. Also in
\cite{AMR} it is
 shown that any cubed $3$--manifold which has all edges of even degree is
geometric
 in the sense of Thurston.

 This research was supported by the Australian Research Council.

\section{Preliminaries}\label{S:Prelim}

In this section we introduce the basic concepts and definitions
 needed for the paper. All $3$--manifolds will be compact and connected.
 All surfaces will be compact but not necessarily connected nor orientable.
All maps
 will be assumed smooth.

\begin{defn}

A properly immersed compact surface not equal to a $2$--sphere, projective plane
 or disk is {\it incompressible} if the induced map of fundamental groups
of the surface
 to the $3$--manifold is injective. We say that the surface is {\it boundary
incompressible}
if no essential arc on the surface with ends in its boundary is homotopic
keeping its ends
 fixed into the boundary of the $3$--manifold.
\end{defn}

 By the work of Schoen and Yau \cite{SY}, any
incompressible surface can be homotoped to be least area in its homotopy class,
 assuming that some Riemannian metric is chosen on the $3$--manifold so that
the boundary
 is {\it convex} (inward pointing mean curvature). Then by Freedman, Hass,
Scott \cite{FHS},
 after possibly a small perturbation, the least area map becomes a
self-transverse
 immersion which lifts to a collection of embedded planes in the universal
cover of
the $3$--manifold. Moreover any two of these planes which cross, meet in proper
 lines only, so that there are no simple closed curves of intersection.
 We will always assume that an incompressible surface is homotoped
 to satisfy these conditions.

Alternatively, the method of PL minimal surfaces can be used instead of
 smooth minimal surfaces (cf \cite{JR}).
 However it is interesting to arrive at
 this conclusion from other considerations, which we will do in
Section \ref{S:Topo}, for the special
surfaces associated with cubings of non-positive curvature.
 We will denote by {$\cal P$} the collection of planes covering a given
 incompressible surface $F$ in $M$.

\begin{defn}
$F$ satisfies the {\it $k$--plane property} for some positive integer $k$
if any subfamily of $k$ planes in {$\cal P$} has a disjoint pair.
$F$ satisfies the {\it $1$--line property}
 if any two intersecting planes in {$\cal P$} meet in a single line.
Finally, the {\it triple point property} is defined for surfaces $F$ which
already
obey the $1$--line property. This condition states that for any three planes
$P$, $P'$
 and $P''$ of {$\cal P$} which mutually intersect in lines, the three
lines cross in an odd
 (and hence finite) number of triple points.
\end{defn}

\begin{rem}

The triple point property rules out the possibility
that the common stabiliser of the three intersecting planes
is non-trivial. Indeed, the cases of either three disjoint
lines of intersection or three lines meeting in infinitely many
triple points are not allowed. A non-trivial common stabiliser
would require one of these two possibilities.

\end{rem}

\begin{defn}

A $3$--manifold is said to be {\it irreducible} if every embedded $2$--sphere bounds a
$3$--cell.
It is {\it $P^2$--irreducible} if it is irreducible and there are no embedded
$2$--sided projective planes.

\end{defn}

  From now on we will be dealing with $3$--manifolds which are
$P^2$--irreducible.

\begin{defn}

A $3$--manifold will be called {\it Haken} if it is $P^2$--irreducible and
either has non-empty
incompressible boundary or is closed and admits a closed embedded $2$--sided
incompressible surface.

\end{defn}

\begin{defn}

A compact $3$--manifold admits a {\it cubing of non-positive curvature}
(or just a {\it cubing} ) if it can be formed by gluing together a finite
collection of standard Euclidean cubes with the following conditions:
\begin{itemize}
\item[-] Each edge of the resulting cell structure has degree at least four.

\item[-] Each link of a vertex of the cell structure is a
triangulated $2$--sphere so that any loop of edges of length three
 bounds a triangle, representing the intersection with a single cube,
 and there are no cycles of length less than three, other than a single
edge traversed twice
with opposite orientations.
\end{itemize}
\end{defn}

\begin{rem}

Note the conditions on the cubing are just a special case of Gromov's
CAT$(0)$ conditions.

\end{rem}

\begin{defn}

The {\it canonical surface} $S$ in a cubed $3$--manifold is
formed by gluing together the finite collection of
squares, each obtained as the intersection of a plane with a cube, where
the plane is midway between a pair
of parallel faces of the cube.

\end{defn}

\begin{rems}\begin{enumerate}\item  
We consider cubed $3$--manifolds which are closed or with boundaries
consisting of
incompressible tori and Klein bottles. In the latter case, the boundary
surfaces
are covered by square faces of the cubes. Such cubings are very useful for
studying
knot and link complements.

\item  In \cite{AR1} it is sketched why the canonical
surface is incompressible and satisfies the $4$--plane, $1$--line and triple
point properties.
 We prove this in detail in the next section for completeness.

\item Regions complementary to $S$ are cones on unions of squares obtained
by canonically decomposing each triangle in the link of a vertex into $3$
squares.

\item The complementary regions of $S$ in (3) are polyhedral cells. Each
such polyhedron
$\Pi$ has boundary determined by a graph $\Gamma_\Pi$ whose edges
correspond to arcs
of double points of $S$ and whose vertices are triple points of $S$. By
construction, there is a unique vertex $v$ in the original cubing which
lies in the
centre of $\Pi$: the graph $\Gamma_\Pi$ is merely the graph on $S^2$ dual
to the
triangulation determined by the link of $v$ in the cubing. The conditions
on the
cubing translate directly into the following statements concerning the
graph $\Gamma_\Pi$.

\begin{itemize}\item[(a)] Every face has degree at least $4$.

\item[(b)] Every embedded loop on $S^2$ meeting $\Gamma_\Pi$ transversely
in exactly two points bounds a disk
cutting off a single arc of $\Gamma_\Pi$.

\item[(c)]  Every embedded loop on $S^2$ meeting $\Gamma_\Pi$ transversely
in exactly three points bounds a disk
which contains a single vertex of degree $3$ of $\Gamma_\Pi$. So the part
of $\Gamma_\Pi$ inside
this disk can be described as a `Y'.
\end{itemize}\end{enumerate}
\end{rems}

\begin{nota}

We will use $S$ to denote a subset of $M$, ie, the image of the canonical
surface and ${\bar S}$
to denote the domain of a map $f\co {\bar S} \rightarrow M$ which has image
$f({\bar S}) = S$.

A similar notational convention applies for other surfaces in $M$.

\end{nota}

\section{The canonical surface}\label{S:Canon}

 In this section we verify the properties listed in the final remark in the
previous section.

\begin{thm}\label{T:1}

 The canonical surface $S$ in a cubed $3$--manifold $M$ is incompressible
and satisfies the $4$--plane, $1$--line and triple point properties.
 Moreover $S$ is covered by a collection of embedded planes in the
universal covering of $M$
 and two such planes meet at most in a single line. Also two such lines meet
in at most a single point.

\end{thm}

\begin{proof}
We show first that $S$ is incompressible.  Of course this follows by
standard techniques, by thinking of $M$ as having a polyhedral metric
of non-positive curvature and using the Cartan--Hadamard Theorem to
identify the universal
 covering with $\bf R^3$ (cf \cite{BGS}). Since $S$ is totally geodesic and
geodesics diverge
 in the universal covering space, we see that $S$ is covered by a collection
 of embedded planes {$\cal P$}.

However we want to use a direct combinatorial argument which generalises to
 situations in the next section where no such metric is obvious on $M$. Suppose
 that there is an immersed disk $D$ with boundary $C$ on $S$.  Assume that
$D$ is in
general position relative to $S$, so that the inverse image of $S$ is a
collection $G$
 of arcs and loops in the domain of the map ${\bar D}\rightarrow M$ with
image $D$.
$G$ can be viewed as a graph with vertices of degree four in the interior
of ${\bar D}$.
 Let $v$ be the number of vertices, $e$ the number of edges and $f$ the
number of faces
of the graph $G$, where the faces are the complementary regions of $G$ in
{$\bar D$}. We assume
 initially that these regions are all disks.

An Euler characteristic argument gives that $v - e + f =1$ and so since $2v
\le e$,
there must be some faces
with degree less than four.  We define some basic homotopies on
 the disk $D$ which change $G$ to eventually decrease the number of
vertices or edges.
First of all assume there is a region in the complement of $G$ adjacent to
$C$ with two or three
 vertices. In the former case we have a $2$--gon $D'$ of {$\bar D$} with one
boundary arc on $C$
 and the other on $G$. So $D'$ has interior disjoint from $S$
 and its boundary lies on $S$. But by definition of $S$, any such a $2$--gon
can be
homotoped into a double arc of $S$. For the $2$--gon is contained in a cell
in the closure of the complement of $S$.
The cell has a polyhedral structure which can be described as the cone on
the dual cell
 decomposition of a link of a vertex of the cubing.
 The two arcs of the $2$--gon can be deformed into the $1$--skeleton of the
link and then define
a cycle of length two. By definition such a cycle is an edge taken twice
in opposite directions. We now homotop $D$ until $D'$ is pushed into the
double arc of $S$
 and then push $D$ slightly off this double arc. The effect is to eliminate
the $2$--gon $D'$,
 ie one arc or edge of $G$ is removed.

 Next assume there is a region $D''$ of the complement of $G$ bounded by
three arcs,
two of which are edges of $G$ and one is in $C$. The argument is very
similar to that in
the previous paragraph. Note that when the boundary of $D$ is pushed into
the $1$--skeleton of
the link of some vertex of the cubing then it gives a $3$--cycle which is the
boundary of a triangle representing the intersection of the link with a
single cube.
Therefore we can slide $D$ so that $D''$ is pushed into the triple point of
$S$ lying at the
centre of this cube. Again by perturbing $D$ slightly off $S$, $D''$ is
removed and $G$
 has two fewer edges and one fewer vertex.

   Finally to complete the argument we need to discuss how to deal with
internal
 regions which are $1$--gons, $2$--gons or $3$--gons. Now
$1$--gons cannot occur, as there are no $1$--cycles in the link of a vertex.
$2$--gons can be
eliminated as above. The same move as described above on $3$--gons has the
effect of
inverting them, ie, moving one of the edges of the $3$--gon past the
opposite vertex.
This is enough to finish the argument by the following observations.

 First of all consider an arc of $G$ which is a path with both ends on $C$
and passes through
vertices by using opposite edges at each vertex (degree $4$),
ie, corresponds to a double curve of $D \cap S$.
If such an arc has a self intersection, it is easy to see there are
embedded $2$--gons or $1$--gons
 consisting of subarcs. Choosing an innermost such a $2$--gon or $1$--gon,
then there must be `boundary'
$3$--gons (relative to the subdisk defined by the $2$--gon or $1$--gon)
if there are intersections with arcs. Now push all $3$--gons off such a
$1$--gon or $2$--gon,
starting with the boundary ones. Then we arrive at an innermost $1$--gon or
$2$--gon with no
arcs crossing it and can use the previous method to obtain a contradiction
or to
decrease the complexity of $G$. Similarly
if two such arcs meet at more than one point, we can remove $3$--gons from
the interior of an
innermost resulting $2$--gon and again simplify $G$. Finally if such arcs
meet at one point, we
get $3$--gons which can be made innermost. So in all cases $G$ can be
reduced down to the empty
 graph, with $C$ then lying in a face of $S$. So $C$ is contractible on $S$.

It remains to discuss the situation when some regions in the complement of
$G$ are not
 disks. In this case, there are components of $G$ in the interior of {$\bar
D$}. We simplify such
 an innermost component $G'$ first by the same method as above, working
with the subdisk
{$\bar D^*$} consisting of all the disks in the complement of $G'$ and
separated
from $C$ by $G'$. So we can get rid of $G'$
 and continue
 on until finally all of $G$ is removed by a homotopy of $D$. Note that
once $D$ has no interior
 intersections
 with $S$ then $D$ can be homotoped into $S$ as it lies in a single cell,
which
has the polyhedral structure of the cone on the dual of a
link of a vertex of the cubing. This completes the argument showing that
$S$ is incompressible.

Next we wish to show why $S$ has the $4$--plane, $1$--line and triple point
properties.
 Before discussing this, it is helpful to discuss why
the lifts of $S$ to the universal covering are planes, without using the
polyhedral metric.
Suppose some lift of $S$ to the universal covering {$\cal M$} was not
embedded. We know such a lift
 $P$ is an immersed plane by the previous argument that $S$ is
incompressible. It is easy to
see we can find an immersed disk $D$ with boundary $C$ on $P$ which
represents a $1$--gon. There
 is one vertex where $C$ crosses a double curve of $P$. But the same
argument as in the previous
 paragraph applies to simplify the intersections of $D$ with all the lifts
$P$ of $S$. We get a
contradiction, as there cannot be a $1$--gon with interior mapped into the
complement of $S$.
This establishes that all the lifts $P$ are embedded as claimed.

It is straightforward now to show that any pair of such planes $P,
P^\prime$ which intersect,
 meet in a single line. For if there
is a simple closed curve of intersection, again
the disk $D$ bounded by this curve on say $P$ can be homotoped relative to
the other planes
to get a contradiction. Similarly if there are at least two lines
of intersection of $P$ and $P^\prime$ then there is a $2$--gon $D$ with
boundary arcs on $P$
 and $P^\prime$.
Again we can deform $D$ to push its interior off {$\cal P$} giving a
contradiction.
This establishes the $1$--line property.

The $4$--plane and triple point properties follow once we can show that any
three planes of ${\cal P}$ which mutually intersect, meet in a single triple point.
For then if four planes all met in pairs, then on one of the planes we
would see
three lines all meeting in pairs. But this implies there is a $3$--gon
between the lines and the same disk simplification argument as above, shows
that this is impossible.
 There are two situations which could occur for three mutually crossing planes
 $P$, $P'$ and $P''$. First of all, there could be no triple points at all
between the
 three planes. In this case the $3$--gon $D$ with three boundary arcs
joining the three
lines of intersection on each of the planes can be used to give a
contradiction. This follows
by the same simplification argument, since the $3$--gon can be homotoped to
have interior
disjoint from ${\cal P}$.
Secondly there could be more than one triple point between the planes. But
in this
 case, in say $P$, we would see two lines meeting more than once. Hence
there would be
$2$--gons $D$ in $P$ between these lines. The interiors of such $2$--gons can
be homotoped
off ${\cal P}$ and the resulting contradiction
completes the argument for all the
properties claimed in the theorem.\end{proof}

\begin{rems}\begin{enumerate}
\item In the next section we will define a class of $3$--manifolds
which are almost
cubed. These do not have nice polyhedral metrics arising from a simple
construction like
cubings but the same methods will work as in the above theorem. This is
the basis of what we call a combinatorial version of Dehn's lemma and the
loop theorem.

\item In \cite{ARS}, other generalisations of cubings are given, where
the manifold behaves as if `on average' it has non-positive curvature.
Again the technique of the above theorem applies in this situation, to
deduce incompressibility of particular surfaces.

\item A key factor in making the above method work is
that there must always be faces of the graph $G$ which are $2$--gons or
$3$--gons.
In particular the Euler characteristic argument to show existence of such
regions breaks down once $D$ is a $4$--gon!
\end{enumerate}\end{rems}

\begin{defn}

An immersed surface $S$ is called {\it filling} if the closures of the
complementary regions of $S$
in $M$ are all cells, for all least area maps in the homotopy class of $S$,
for any metric on $M$.

\end{defn}

  It is trivial to see that for the canonical surface $S$ in a cubed
$3$--manifold $M$,
all the closures of the complementary regions are cells.
 A little more work checks that $S$ is actually filling. In fact, one way
to do this
is to observe that any essential loop in $M$ can be homotoped to a geodesic
$C$ in the
polyhedral metric defined by the cubing.
Then this geodesic lifts to a line in the universal covering {$\cal M$}. A
geodesic line will meet
 the planes over $S$ in single points or lie in such a plane.
Note that the lines of intersection between the planes are also geodesics.
So if the
given geodesic line lies in some plane,
then by the filling property, it meets some other line of intersection in a
single
point. Hence it meets the corresponding plane in one point.

 A homotopy of $S$ to a least
area surface $S^\prime$ relative to some metric, will lift to a proper
homotopy between
collections of embedded planes ${\cal P}$ and ${\cal P^\prime}$
 in the universal covering. This latter homotopy cannot remove such
essential intersections
between the given geodesic line and some plane (all points only move a
bounded distance,
whereas the ends of the line are an unbounded distance from the plane and
on either side of it).
So we conclude that any essential loop must intersect $S$.

Therefore a complementary domain to $S^\prime$ must have fundamental group
with trivial
image in  ${\pi_1}(M)$.
The argument of \cite{HRS} shows that for a least area map of an incompressible
surface, all the
complementary regions are $\pi_1$--injective. So such regions must be cells,
as a cubed
manifold is irreducible.

\begin{rem}

Note that $P^2$--irreducibility for a cubed manifold can also be shown directly,
since we can apply the same argument as
above to simplify the intersections of an immersed sphere or projective plane
 with the canonical surface and
eventually by a homotopy, achieve that the sphere or projective plane lies
in a complementary cell.

\end{rem}

  In \cite{AR1} it is observed that a $3$--manifold has a cubing of non-positive
curvature
if and only if it has a filling incompressible surface satisfying the
$4$--plane, $1$--line and
 triple point properties. This follows immediately from the work of Hass
and Scott \cite{HS}

\section{Combinatorial Dehn's lemma and the loop theorem}\label{S:Comb}

 Our aim here is to define a class of $3$--manifolds which are almost
 cubed and for which one can still verify that the canonical surface
is imcompressible and satisfies similar properties to that for cubings.
 In fact the canonical surface here satisfies the $4$--plane, $1$--line and
triple
point properties, but not necessarily the filling property. So as in
\cite{HRS},
 the closures of the complementary regions of $S$ are $\pi_1$--injective
handlebodies.
These surfaces naturally arise in \cite{AR3}, where we investigate surgeries on
certain
 classes of simple alternating links containing such closed surfaces in
their complements.

\begin{defn}

Suppose that $S$ is an immersed closed surface in a compact
3--manifold $M$. We say that the closure of a connected component of the
complement in $S$
of the double curves of $S$, is a {\it face} of $S$.

\end{defn}

\begin{defn}

Suppose that $D$ is an immersed disk in $M$ with boundary $C$ on an immersed
closed surface $S$ and interior of $D$ disjoint from $S$.
Assume also that there are no such disks $D$ in $M$ for which $C$ crosses the
double curves of $S$ exactly once. We say that {\it $D$ is homotopically
trivial
relative to $S$}, if one of the following three situations hold:
\begin{enumerate}
\item If $C$ has no intersections with the double curves of $S$,
then $D$ can be homotoped into a face of $S$, keeping its boundary fixed.

\item Any $2$--gon $D$ (ie, $C$ meets the singular set of $S$ in two
points) is homotopic into a double curve of $S$, without changing the
number of intersections
of $C$ and the double curves of $S$ and keeping $C$ on $S$.

\item Any $3$--gon $D$ is homotopic into a triple point of $S$, without
changing the number of intersections
of $C$ and the double curves of $S$ and keeping $C$ on $S$.
\end{enumerate}
\end{defn}

\begin{rem}

Note that these conditions occur in Johannson's definition of boundary
patterns in \cite{Jo1}, \cite{Jo2}.

\end{rem}

\begin{thm}\label{T:2}

Assume that $S$ is an immersed closed surface in a compact
3--manifold $M$. Suppose that all the faces of $S$ are polygons with at
least four sides.
Also assume that any embedded disk $D$ with boundary $C$ on $S$ and
interior disjoint from $S$,
with $C$ meeting the double curves two or three times, is homotopically
trivial relative to $S$
and there is no such disk with $C$ crossing the double curves once.

Then $S$ is incompressible with the $4$--plane, $1$--line and triple point
properties.
Moreover $S$ lifts to a collection of embedded planes in the universal
cover of $M$
and each pair of these planes meets in at most a single line. If three planes
mutually intersect in pairs, then they share a single triple point.  Also the
closures of components of the complement of $S$ are $\pi_1$--injective.
Finally if $S$ is incompressible with the $4$--plane, $1$--line and triple
point properties
then $S$ can be homotoped to satisfy the above set of conditions.

\end{thm}

\begin{proof}

  The proof is extremely similar to that for Theorem \ref{T:1} so we only
remark on the ideas. First of all the conditions in the statement of
Theorem \ref{T:2} play the same role as the link conditions in the
definition of a cubing
 of non-positive curvature. So we can homotop disks which have boundary on
$S$ to
 reduce the graph of intersection of the interior of the disk with $S$. In
this way,
 $2$--gons and $3$--gons can be eliminated, as well as compressing disks for
$S$. This
 is the key idea and the rest of the argument is entirely parallel to
Theorem \ref{T:1}.
Note that the closures of the complementary regions of $S$ are
$\pi_1$--injective,
by essentially the same proof as in \cite{HRS}.

The only thing that needs to be carefully checked, is why it suffices to
assume
that only embedded disks in the complementary regions need to be examined, to
see that any possibly singular $n$--gons, for $2 \le n \le 3$, are
homotopically trivial
and there are no singular $1$--gons.

Suppose that we have a properly immersed disk $D$ in a complementary
region, with boundary meeting
the set of double curves of $S$ at most three times and $D$ is not
homotopic into a
 double arc, a triple point or a face of $S$.  If this disk is not
homotopic into
the boundary of the complementary region, we can apply Dehn's lemma and the
loop theorem to
 replace the singular disk by an embedded one. Moreover since the boundary
of the
 new disk is obtained by cutting-and-pasting of the old boundary curve, we
see the
 new curve also meets the set of double curves of $S$ at most three times.
So this
case is easy to handle: it does not happen.

Next assume that the singular disk is homotopic into the boundary surface
$T$ of
the complementary region. (Note we include the possibility here that the
complementary region is a ball
 and $T$ is a $2$--sphere). Let $C$ be the boundary curve of the singular
disk and let $N$
 be a small regular neighbourhood of $C$ in $T$. Thus $C$ is null homotopic
in $T$.
 Notice that there are at most three double
 arcs of $S$ crossing $N$. Now fill in the disks $D'$ bounded by any
contractible boundary
 component $C'$ of $N$ in $T$, to enlarge $N$ to $N'$. Since $C$ shrinks in
$T$, it is easy to see
 by Van Kampen's theorem, that $C$ contracts also in $N'$. Also if $C'$
meets the double
arcs of $S$, we see the picture in $D'$ must be either a single arc or
three arcs meeting
at a single triple point, or else we have found an embedded disk
contradicting our
assumption. For we only need to check that $C'$ cannot meet the double arcs
in at
least four points. If $C'$ did have four or more intersection points with
the singular set of $S$, then one of the double arcs crossing $N$ has both
 ends on $C'$. But this is impossible, as there would be a cycle in the
graph of the
double arcs on $T$, which met the contractible curve $C$ once.

Finally we notice that there must be some disks $D'$ which meet the double
arcs; in fact at
least one point on the end of each double arc in the boundary of $N$ must
be in such a disk.
 For otherwise it is impossible for $C$ to shrink in $N'$, as there is an
essential intersection
at one point with such an arc.
(This immediately shows the possibility that $C$ crossed the double curves
once cannot happen).
 So there are either one or two disks $D'$ with a single arc and
 at most one such disk with three arcs meeting in a triple point. But the
latter case means that
$C$ can be shrunk into the triple point and the former means $C$ can be
homotoped into the double
 arc of $S$ in $N'$ by an easy examination of the possibilities.

 Hence this shows that it suffices to consider only embedded disks
 when requiring the properties in Theorem \ref{T:2}. This is very useful
in applications in \cite{AR3}.

To show the converse, assume we have an incompressible surface which has
the $4$--plane, $1$--line and
triple point properties. Notice that in the paper of Hass and Scott \cite{HS},
the triple point
 property is enough to show that once the number of triple points has been
minimised
for a least area representative of $S$, then the combinatorics of the
surface are rigid.
So we get that $S$ has exactly the properties as in Theorem  \ref{T:2}.
\end{proof}

\begin{rem}

Theorem  \ref{T:2} can be viewed as a singular version of Dehn's lemma
and the loop theorem.
 For we have started with an assumption that there are no embedded disks of
a special
type with boundary on the singular surface $S$ and have concluded that $S$
is incompressible,
ie $\pi_1$--injective. In \cite{ARS} other variants on this theme are
given.
\end{rem}

\begin{defn}

 We say that $M$ is {\it almost cubed} if it is  $P^2$--irreducible and
contains a surface $S$ as in Theorem  \ref{T:2}.

\end{defn}

It is interesting to speculate as to how large is the class of almost cubed
$3$--manifolds.
We do not know of any specific examples of compact $P^2$--irreducible
$3$--manifolds
with infinite $\pi_1$
which are not almost cubed.

\begin{cor}

Suppose that $M$ is a compact $P^2$--irreducible $3$--manifold with boundary,
which is almost cubed,
ie, there is a canonical surface $S$ in the interior of the manifold.
Assume also that the
complementary regions of $S$ include collars of all components of the
boundary. Suppose that a
 handlebody is glued onto each boundary component of $M$ to give a new
manifold $M'$.
If the boundary of every meridian disk, when projected onto $S$,
meets the double curves at least four times,
then $M'$ is almost cubed.

\end{cor}

\begin{proof}
This follows immediately from Theorem  \ref{T:2}, by observing
that since the boundary of every meridian disk meets the double curves at
 least four times, there are no non-trivial $n$--gons in the complement of
 $S$ in $M'$ for $n= 2,3$ and no $1$--gons.  Hence $M'$ is almost cubed, as
$S$ in $M'$ has similar
 properties to $S$ in $M$.\end{proof}

\begin{rem}

Examples satisfying the conditions of
the corollary are given in \cite{AR3}. In particular
 such examples occur for many classes of simple alternating link complements.
In \cite{ALR}, the class of well-balanced alternating links are shown to be
almost
cubed and so the corollary applies.

\end{rem}

\section{Hierarchies}\label{S:Hier}

Our aim in this section is to give a similar treatment of hierarchies, to
that of cubings
and almost cubings.
The definition below is motivated by the special hierarchies used by
Waldhausen
in his solution of the word problem in the class of Haken $3$--manifolds
\cite{Wa1}. Such hierarchies
were extensively studied by Johannson in his work on the characteristic
variety theorem \cite{Jo1}
and also in \cite{Jo2}.

\begin{defn}

A {\it hierarchy} is a collection ${\cal S}$ of embedded compact $2$--sided
surfaces
 ${S_1, S_2,\ldots ,S_k}$,
which are not $2$--spheres, in a compact
 $P^2$--irreducible
$3$--manifold, with the following properties:

\begin{enumerate}
  \item Each $S_i$ has boundary on the union of the previous $S_j$ for $j
\le i-1$.

  \item If an embedded polygonal disk $D$ intersects ${\cal S}$
in its boundary loop $C$ only and $C$ meets the boundary curves of ${\cal S}$
in at most $3$ points, then $D$ is homotopic into either an arc of a
boundary curve,
a vertex or a surface of ${\cal S}$, where $\partial D$ is mapped into
${\cal S}$
throughout the homotopy.

		 \item Assume an embedded polygonal disk $D$ intersects
${\cal S}$
in its boundary loop $C$ only and $C$ has only one boundary arc $\lambda$
on the surface $S_j$,
where $j$ is the largest value for surfaces of ${\cal S}$ met by $C$. Then
$\lambda$
is homotopic into the boundary of $S_j$ keeping the boundary points of
$\lambda$ fixed.
\end{enumerate}
\end{defn}

\begin{rems}\begin{enumerate}\item 
Note that Waldhausen shows that for a Haken $3$--manifold,
 one can always change a given hierarchy into one satisfying these conditions
 by the following simple procedures:

\begin{itemize}\item[-] Assuming that all $S_j$ have been picked for $j<i$, then first arrange
that after
$M$ is cut open along all the $S_j$ to give $M_{i-1}$, the boundary of
$S_i$ is chosen so
that there are no triangular regions cut off between the `boundary pattern'
of $M_{i-1}$
(ie, all the boundary curves of surfaces $S_j$ with $j<i$) and the
boundary curves of $S_i$. This is done by minimising the intersection between
the boundary pattern of $M_{i-1}$ and $\partial S_i$.

 \item[-] It is simple to arrange that $S_i$ is boundary incompressible in
$M_{i-1}$, by performing
 boundary compressions if necessary. So there are no $2$--gons between $S_i$
and $S_j$ for
any $j<i$.

 \item[-] Finally one can see that there are no essential triangular disks in
$M_{i-1}$, with one
 boundary arc on $S_i$ and the other two arcs on surfaces $S_j$ for values
$j \le i-1$,
 by the boundary incompressibility of $S_i$ as in the previous step.
\end{itemize} 

 \item  Notice we are not assuming the hierarchy is complete, in the sense
that the complementary
 regions are cells (in the case that $M$ is closed) or cells and collars of
the boundary
(if M is compact with incompressible boundary).

 \item  The next result is a converse statement, showing that the
conditions above imply that
the surfaces do form a hierarchy.
\end{enumerate}
\end{rems}

\begin{thm}\label{T:3}

Assume that ${S_1, S_2,\ldots , S_k}$ is a sequence ${\cal S}$ of embedded
compact 2--sided surfaces,
none of which are $2$--spheres, in a compact $P^2$--irreducible
$3$--manifold $M$ with the following properties:
\begin{enumerate}
   \item Each $S_i$ has boundary on the union of the previous $S_j$ for $j
\le i-1$.

  \item If an embedded polygonal disk $D$ intersects ${\cal S}$
in its boundary loop $C$ only and $C$ meets the boundary curves of ${\cal S}$
in at most $3$ points, then $D$ is homotopic into either an arc of a
boundary curve,
a vertex or a surface of ${\cal S}$.

		\item Assume an embedded polygonal disk $D$ intersects
${\cal S}$
in its boundary loop $C$ only and $C$ has only one boundary arc $\lambda$
on the surface $S_j$,
where $j$ is the largest value for surfaces of ${\cal S}$ met by $C$. Then
$\lambda$
is homotopic into the boundary of $S_j$ keeping the boundary points of
$\lambda$ fixed.
\end{enumerate}

Then each of the surfaces $S_i$ is incompressible and boundary
incompressible in the cut open manifold $M_{i-1}$ and ${\cal S}$
forms a hierarchy for $M$ as above.

\end{thm}

\begin{proof}

The argument is very similar to those for Theorems  \ref{T:1} and
\ref{T:2} above, so we outline the
modifications needed.

Suppose there is a compressing or boundary compressing disk $D$ for one of
the surfaces $S_i$.
 We may assume that all the previous $S_j$ are incompressible and boundary
incompressible by
 induction. Consider $G$, the graph of intersection of $D$ with the
previous $S_j$,
pulled back
to the domain of $D$. Then $G$ is a degree three graph; however at each
vertex there is a `$T$' pattern
 as one of the incident edges lies in some $S_j$ and the other two in the
same surface
 $S_k$ for some $k<j$.
We argue
 that the graph $G$ can be simplified by moves similar to the ones in
Theorems  \ref{T:1} and  \ref{T:2}.

First of all, note that by an innermost region argument, there must be
either an innermost
$0$--gon, $2$--gon or a
triangular component of the closure of the complement of $G$.  For we can
cut up $D$ first using
the arcs of intersection with $S_1$, then $S_2$ etc. Using the first
collection of arcs, there is
clearly an outermost $2$--gon region in $D$. Next the second collection of
arcs is either disjoint
 from this $2$--gon or there is an outermost $3$--gon. At each stage, there
must always
be an outermost $2$--gon or $3$--gon. (Of course any simple closed curves of
intersection just
 start smaller disks which can be followed by the same method. If such a
loop is isolated, one gets
 an innermost $0$--gon which is readily eliminated, by assumption).

 By supposition, such a $2$--gon or $3$--gon can be homotoped into either a
boundary curve or
 into a boundary vertex of some $S_j$. We follow this by slightly deforming
the map
to regain general position. The complexity of the graph is defined
 by listing lexicographically the numbers of vertices with a particular
label. The
label of each vertex is
given by the subscript of the first surface of $\cal S$ containing the vertex.
(cf \cite{Jo2} for a good discussion of
this lexicographic complexity). The homotopy above can be readily seen to
reduce the complexity of the graph.
 Note that the hypotheses only refer to embedded
 $n$--gons but as in the proof of Theorem  \ref{T:2} it is easy to show
that if there are only
trivial embedded $n$--gons for $n<4$, then the same is true for immersed
$n$--gons, using Dehn's
lemma and the loop theorem. Similarly, the hypothesis (3) of the
theorem
can be converted to a statement about embedded disks, using Dehn's lemma
and the
loop theorem, since cutting open the manifold using the previous surfaces,
converts the polygon into a $2$--gon. This completes the proof of
Theorem  \ref{T:3}.\end{proof}

\begin{exam}

Consider the Borromean rings complement $M$ in the $3$--sphere. We can find
such a hierarchy
easily as follows:

Start with a peripheral torus as $S_1$, ie, one of the boundary tori of
the complement $M$.
Next choose $S_2$ as an essential embedded disk with boundary on $S_1$,
with a tube attached to
 avoid the intersections with one of the other components $C_1$ of $B$, the
Borromean rings.
Now cut $M$ open along $S_1$ and $S_2$ to form a collar of $S_1$ and
another component $M_2$.
 It is easy to see that $M_2$
 is a genus $2$
 handlebody and $B$ has two components $C_1$ and $C_2$ inside this, forming
a sublink
$B'$.  Moreover
 these two loops are generators of the fundamental group of $M_2$, since
they are dual
(intersect in single points) to two disjoint meridian disks for $M_2$.
Finally the
loops are readily seen to be linked once in $M_2$.

Now we can cut $M_2$ open using a separating meridian disk with a handle
attached, so
 as to avoid $B'$. This surface $S_3$ can also be viewed as another
spanning disk with
boundary on $S_1$ having a tube attached to miss the other component $C_2$
of $B'$. So
 there is a nice symmetry between $S_2$ and $S_3$. Finally we observe that
when $M$ is cut
 open along $S_1$, $S_2$ and $S_3$ to form $M_3$, this is a pair of genus
two handlebodies, each of
 which contains one of $C_1$ and $C_2$ as a core curve for one of the handles,
plus the collar of $S_1$.  So we can
 choose $S_4$ and $S_5$ to be non-separating meridian disks for these
handlebodies, disjoint
 from $B$, so that $M_5$ consists of collars of the three boundary tori.

  Notice that  the `boundary patterns' on each of these tori, ie, the
boundary curves of
 the surfaces in the hierarchies, consist of two contractible simple closed
curves and
four arcs, two of which have boundary on each of these loops. The pairs of
arcs are `parallel',
 in that the whole boundary pattern divides each torus into six regions,
two $0$--gon disks,
two $4$--gon disks and two annuli (cf Figure 4 in \cite{AR2}).

As a corollary then to Theorem \ref{T:3}, we observe that any non-trivial
surgery on each of the
two components $C_2$ and $C_3$ of $B$ gives meridian disks $S_6$ and $S_7$
which
 meet the boundary pattern
 at least four times. Hence these surfaces form a hierarchy for the
surgered manifold and so all
 such surgeries on $C_2$ and $C_3$ have the result
 that the peripheral torus $S_1$ remains incompressible.

\end{exam}

\begin{rems}\begin{enumerate}\item
This can be proved by other methods but the above argument is
particularly revealing,
 not using any hyperbolic geometry.

\item A similar example of the Whitehead link is discussed in \cite{AR2} and
it is interesting
to note the boundary pattern found there (Figure 4) is exactly the same as here.

\item A very significant problem is to try to use the above
characterisation of hierarchies
 to find some type of polyhedral metric of non-positive curvature, similar
to cubings.
This would give a polyhedral approach to Thurston's uniformisation Theorem
for Haken
manifolds \cite{Th}.
\end{enumerate}
\end{rems}

\section{Topological rigidity of cubed manifolds}\label{S:Topo}

  In this section we give a different approach to the result of Hass and
Scott \cite{HS}
that if a closed $P^2$--irreducible $3$--manifold is homotopy equivalent to a
cubed $3$--manifold
 then the manifolds are homeomorphic. This can be viewed as a polyhedral
version of Mostow
 rigidity, which says that complete hyperbolic manifolds of finite volume
which are
 homotopy equivalent are isometric, in dimensions greater than $2$. We
refer to this
 as topological rigidity of cubed $3$--manifolds. Our aim here is to show
that this rigidity
theorem can be shown without resorting to the least area methods of
Freedman Hass Scott
\cite{FHS}, but can be obtained by a direct argument more like Waldhausen's
original proof of
 rigidity for Haken $3$--manifolds \cite{Wa2}. Note that various generalisations
of Hass and Scott's
theorem have been obtained recently by Paterson \cite{P1}, \cite{P2} using
different
methods.


\begin{thm}\label{T:4}
Suppose that $M$ is a compact $3$--manifold with a cubing of non-positive
curvature
and $M'$ is a closed $P^2$--irreducible $3$--manifold which is homotopy
equivalent to $M$.
Then $M$ and $M'$ are homeomorphic.
\end{thm}

\begin{proof}
Note that the cases where either $M$ has incompressible boundary or
is non-orientable,
are not so interesting, as then $M$ and $M'$ are Haken and the result
follows by Waldhausen's
theorem \cite{Wa2}. So we restrict attention to the case where $M$ and $M'$ are
closed and orientable.

Our method is a mixture of those of \cite{HS} and \cite{Wa2} and we
indicate the
steps, which
are all quite standard techniques.

\noindent{\bf Step 1}\qua
 By Theorem \ref{T:1}, if $S$ is the canonical surface for the cubed
manifold
$M$ and $M_S$
 is the cover corresponding to the fundamental group of $S$, then $S$ lifts
to an embedding
denoted by $S$ again in $M_S$. For by Theorem \ref{T:1}, all the lifts of
$S$ to the universal covering
${\cal M}$ are embedded planes and $\pi_1(S)$ stabilises one of these planes
 with quotient the required
lift of $S$. Let $f\co M' \rightarrow M$ be the homotopy equivalence and
assume that $f$ has been
perturbed
 to be transverse to $S$. Denote the immersed surface $f^{-1}(S$) by $S'$.
Notice that
 $f$ lifts
to a map ${\tilde f}\co  {\cal M'} \rightarrow {\cal M}$ between universal
covers and so all the
 lifts of $S'$ to ${\cal M'}$ are properly
embedded non-compact surfaces. In fact, if $M'_S$ is the induced cover of
$M'$ corresponding
 to $M_S$,
ie, with fundamental group projecting to $f^{-1}(\pi_1(S))$, then there
is an
 embedded lift, denoted
 $S'$ again, of $S'$ to $M'_S$, which is the inverse image of the embedded
lift of $S$.

The first step is to surger $S'$ in $M'_S$ to get a copy of $S$ as the
result. We will be able
to keep some of the nice properties of $S$ by this procedure, especially
the $4$--plane property.
 This will enable us to carry out the remainder of the argument of Hass and
Scott quite
 easily. For convenience, we will suppose that $S$ is orientable.
The non-orientable case is not
 difficult to derive from this; we leave the details to the reader. (All
that is necessary is to
 pass to a $2$--fold cover of $M_S$  and $M'_S$, where $S$ lifts to its
orientable double covering surface.)

Since $f$ is a homotopy equivalence, so is the lifted map $f_S\co M'_S
\rightarrow M_S$. Hence
 if $S'$ is not
 homeomorphic to $S$, then the induced map on fundamental groups of the
inclusion
 of $S'$ into $M'_S$
has kernel in $\pi_1(S')$. So we can compress $S'$ by Dehn's lemma and the
loop
theorem. On the other hand, the ends of $M_S$ pull back to ends of $M'_S$
and a properly
 embedded line going between the ends
of $M_S$ meeting $S$ once, pulls back to a similar line in $M'_S$ for $S'$.
Hence $S'$
 represents a non-trivial
 homology class in $M'_S$ and so cannot be completely compressed. We
conclude that $S'$
 compresses to
 an incompressible surface $S''$ separating the ends of $M'_S$.

 Now we claim that any component $S^*$ of $S''$ which is
homologically non-trivial,
 must be homeomorphic
 to $S$. Also the inclusion of $S^*$ induces an isomorphism on fundamental
groups to $M'_S$.
 The argument
 is in \cite{FHS}, for example, but we repeat it for the benefit of the
reader. The homotopy equivalence
between $M'_S$ and $S$ induces a map $f\co  S^* \rightarrow S$ which is
non-zero on second
homology. So $f$ is homotopic to a finite sheeted covering. Lifting $S^*$ to
the corresponding finite sheeted cover of $M'_S$, we get a number of copies
of $S^*$, if the map
 $S^* \rightarrow M'_S$ is not a homotopy equivalence. Now the
different lifts of $S^*$
 must all separate the two ends of the covering of $M'_S$ and so are all
homologous. (Note as the
second homology is cyclic, there are exactly two ends). But then
any compact region between these lifts projects onto $M'_S$ and so $M'_S$
is actually compact,
unless $S$ is non-orientable, which has been ruled out.

 Finally to complete this step, we claim that if $S^*$ is projected to $M'$
and then lifted to
 {$\cal M'$},
then the result is a collection of embedded planes ${\cal P^*}$ satisfying
the $4$--plane property.
 Notice
first of all that all the lifts ${\cal P'}$ of $S'$ to ${\cal M'}$ satisfy
the `$4$--surface' property.
In other words,
 if any subcollection of four components of ${\cal P'}$ are chosen,
 then there must be a disjoint pair. This is evident
 as $S'$ is the pull-back of $S$ and so ${\cal P'}$ is the pull-back of
${\cal P}$.
 Then the $4$--plane property clearly
pulls-back to the `$4$--surface' property as required.

Now we claim that as $S'$ is surgered and then a component $S^*$ is chosen,
this can be done
so that the $4$--surface property remains valid. For consider some disk $D$
used to surger the
embedded lift $S'$ in $M'_S$. By projecting to $M'$ and lifting to ${\cal M'}$,
 we have a family of embedded
disks surgering the embedded surfaces ${\cal P'}$.

  It is sufficient to show that one such $D$ can be selected so as to miss
all the surfaces $P''$ in
 ${\cal P'}$ which are disjoint from a given surface $P'$ containing the
boundary of $D$,
as the picture
in the universal covering is invariant under the action of the covering
translation group.
This is similar to the argument in \cite{HRS}. First of all, if $D$ meets any
such a surface $P''$ in a
loop which is non-contractible on $P''$, we can replace $D$ by the subdisk
 bounded by this loop. This
 subdisk has fewer curves of intersection with ${\cal P'}$ than the original.
Of course the subdisk may not be disjoint
from its translates under the stabiliser of $P''$. However we
can fix this up at the end of the argument. We relabel $P''$ by $P'$ if
this step is necessary.

Suppose now that $D$ meets any surface $P''$ disjoint from $P'$ in loops $C$
 which are contractible on $P''$. Choose an innermost such a loop. We would
like to do a
simple disk swap and reduce the number of such surfaces $P''$ disjoint from
$P'$ met by $D$.
Note we do not care if the number of loops of intersection goes up during
this procedure.
 However we must be careful that no new planes are intersected by $D$. So
suppose that $C$ bounds
 a disk $D''$ on $P''$ met by some plane $P_1$ which is disjoint from $P'$,
but $D$
 does not already meet $P_1$. Then clearly $P_1$ must meet $D''$ in a
simple closed
 curve in the interior of $D''$. Now we can use the technique of Hass and
Scott to eliminate
 all such intersections in $P'$. For by the $4$--surface property, either
such simple closed curves
 are isolated (ie not met by other surfaces) or there are disjoint
embedded arcs where
 the curves of intersection of the surfaces cross an innermost such loop.
But then we can
start with an innermost such $2$--gon between such arcs and by simple
$2$--gon moves, push all
the arcs equivariantly outside the loop. At each stage we decrease the
number of
triple points and eventually can eliminate the contractible double curves.
 The conclusion is that eventually we can pull $D$ off all the surfaces
disjoint from $P'$.

  Finally to fix up the disk $D$ relative to the action of the stabiliser
of $P'$, project
$P'$ to the compact surface in $M'_S$. Now we see that $D$ may project to
an immersed disk, but all the lifts of this
immersed disk with boundary on $P'$ are embedded and disjoint from the
surfaces in ${\cal P'}$ which
miss $P'$. We can now apply Dehn's lemma and the loop theorem to replace
the immersed disk
by an embedded one in $M'_S$. This is obtained by cutting and pasting, so
it follows
immediately that any lift of the new disk with boundary on $P'$ misses all
the
surfaces which are disjoint from $P'$ as desired. This completes the first
step of the argument.

\noindent{\bf Step 2}\qua
The remainder of the argument follows that of Hass and Scott closely. By
Step 1 we have a
component $S^*$ of the surgered surface which gives a subgroup of
$\pi_1(M')$ mapped by $f$
isomorphically
to the  subgroup of $\pi_1(M)$ corresponding to $\pi_1(S)$. Also $S^*$ is
embedded in the cover $M'_S$
 and all the lifts
 have the $4$--plane property in {$\cal M'$}.
  All that remains is to use Hass and Scott's triple point cancellation
technique to get rid of
 redundant triple points and simple closed curves of intersection between
the planes over $S^*$.
Eventually we get a new surface, again denoted by $S^*$, which is changed
by an isotopy in $M'_S$
 and has
 the $1$--line and triple point properties. It is easy then to conclude that
an equivariant
homeomorphism between {$\cal M'$} and {$\cal M$} can be constructed.

Note that this case is the easy one in Hass and Scott's paper, as the
triple point property
means that triangular prism regions cannot occur and so no crushing of such
regions is necessary.\end{proof}

\begin{cor}

If $f\co M' \rightarrow M$ is a homotopy equivalence and $M$ has an immersed
incompressible surface $S$
satisfying the $4$--plane property, then the method of Theorem \ref{T:4}
shows that
$M'$ also has an immersed incompressible
surface $S'$ satisfying the $4$--plane property and $f$ induces an
isomorphism from
 $\pi_1(S')$ to $\pi_1(S)$.\endproof
\end{cor}

\begin{rem}
This can be shown to be true by least area methods also, but it is
interesting to have
alternative combinatorial arguments. Least area techniques give the result
also for
the $k$--plane property, for all $k$. However there is no direct information
about how
the surface is pulled back, as in the above argument.

\end{rem}

%
%

\Addresses\recd


\begin{thebibliography}



\bibitem{ALR} {\bf I\,R Aitchison}, {\bf E Lumsden}, {\bf J\,H Rubinstein},
{\it Cusp structure of alternating links},
Invent. Math. 109 (1992) 473--494

\bibitem{AMR} {\bf I\,R Aitchison}, {\bf S Matsumoto}, {\bf J\,H Rubinstein},
{\it Immersed surfaces in cubed manifolds},
Asian J. of Math. 1 (1997) 85--95

\bibitem{AR1} {\bf I\,R Aitchison}, {\bf J\,H Rubinstein},
{\it An introduction to  polyhedral metrics  of non-positive curvature},
from: ``Geometry of Low Dimensional Manifolds: 2'',
London Math. Soc. Lecture Notes Series
151, Cambridge University Press (1990) 127--161

\bibitem{AR2} {\bf I\,R Aitchison}, {\bf J\,H Rubinstein},
{\it Incompressible surfaces and the topology of 3--dimensional manifolds},
J. Aust. Math. Soc. (Series A) 55 (1993) 1--22

\bibitem{AR3} {\bf I\,R Aitchison}, {\bf J\,H Rubinstein},
{\it Immersed incompressible surfaces in alternating link complements},
in preparation

\bibitem{ARS}  {\bf I\,R Aitchison}, {\bf J\,H Rubinstein}, {\bf M Sageev},
{\it Separability properties of some immersed incompressible surfaces in
3--manifolds},
in preparation

\bibitem{BGS} {\bf W Ballmann}, {\bf M Gromov}, {\bf V Schroeder},
{\it Manifolds of non-positive curvature},
Progress in Math. 61, Birkhauser (1985)

\bibitem{BH} {\bf S Bleiler}, {\bf C Hodgson},
{\it Spherical space forms and Dehn surgery}, from:
``Knots 90'', de Gruyter (1992) 425--434

\bibitem{FHS} {\bf M Freedman}, {\bf J Hass}, {\bf P Scott},
{\it Least area incompressible surfaces in 3--manifolds},
Invent. Math. 71 (1983) 609--642

\bibitem{HRS} {\bf J Hass}, {\bf J\,H Rubinstein}, {\bf P Scott},
{\it Compactifying covers of closed 3--manifolds},
J. Diff. Geom. 30 (1989) 817--832

\bibitem{HS} {\bf J Hass}, {\bf P Scott},
{\it Homotopy equivalences and homeomorphisms of 3--manifolds},
Topology, 31 (1992) 493--517

\bibitem{Jo1} {\bf K Johannson},
{\it Homotopy equivalences of 3--manifolds with boundaries},
Lecture Notes in Math. 761, Springer (1978)

\bibitem{Jo2} {\bf K Johannson},
{\it On the loop and sphere theorem},
from:
``Low Dimensional Topology'',
Proc. of the Conf. at Univ. of Tennessee,
Conf. Proc. and Lecture Notes on Geom. and Top. 3, International
Press (1994) 47--54

\bibitem{JR} {\bf W Jaco}, {\bf J H Rubinstein},
{\it PL minimal surfaces in 3--manifolds},
J. of Diff. Geom. 27 (1988) 493--524

\bibitem{Mo}  {\bf L Mosher},
{\it Geometry of cubulated 3--manifolds},
Topology, 34 (1995) 789--813

\bibitem{P1} {\bf J Paterson},
{\it The classification of 3--manifolds with
 4--plane quasi-finite immersed surfaces},
preprint

\bibitem{P2} {\bf J Paterson},
{\it Immersed one-to-one projective surfaces in 3--manifolds},
preprint

\bibitem{SY} {\bf R Schoen}, {\bf S\,T Yau},
{\it Existence of incompressible minimal surfaces and the
topology of 3--dimensional manifolds with non-negative scalar curvature},
Annals of Math. (2)
110 (1979) 127--142

\bibitem{Th} {\bf W\,P Thurston},
{\it The geometry and topology of 3--manifolds},
Princeton University (1978)

\bibitem{Wa1} {\bf F Waldhausen},
{\it On irreducible 3--manifolds which are sufficiently large},
Annals of Math. 87 (1968) 56--88

\bibitem{Wa2} {\bf F Waldhausen},
{\it The word problem in fundamental groups
 of sufficiently large 3--manifolds},
Annals of Math. 88 (1968) 272--280

\end{thebibliography}
\end{document}